\documentclass[12pt]{amsart}
\usepackage{amscd,amssymb}
\newtheorem{thm}{Theorem}[section]
\newtheorem{corol}[thm]{Corollary} 
\newtheorem{lemma}[thm]{Lemma}
\newtheorem{prop}[thm]{Proposition}

\newtheorem{question}[thm]{Question}
\theoremstyle{definition}
\newtheorem{defin}[thm]{Definition}

\theoremstyle{remark}
\newtheorem{remark}[thm]{Remark}

\numberwithin{equation}{section}

\def\e{{\varepsilon}}
\def\RUC{{\hbox{RUC\,}^b}}

\def\Is{{\mathrm{Is}}\,}

\def\QED{\nobreak\quad\ifmmode\roman{Q.E.D.}\else{\rm Q.E.D.}\fi}

\def\sS{{\mathcal S}}
\def\sN{{\mathcal N}}
\def\sP{{\mathcal P}}
\def\obr{^{-1}}
\def\sbs{\subset}
\def\a{\alpha}
\def\o{\omega}

\begin{document}

\title[Weakly almost periodic groups]
{A note on the precompactness of weakly almost periodic groups}

\author[M.G. Megrelishvili]{Michael G. Megrelishvili}
\address{Department of Mathematics, 
Bar-Ilan University, 52900 Ramat-Gan, Israel}
\email{megereli@macs.biu.ac.il}
\urladdr{http://www.cs.biu.ac.il/$^\sim$megereli}%
\author[V.G. Pestov]{Vladimir G. Pestov}
\address{School of Mathematical and Computing Sciences,
Victoria University of Wellington, P.O. Box 600, Wellington,
New Zealand}
\email{vova@mcs.vuw.ac.nz}
\urladdr{http://www.mcs.vuw.ac.nz/$^\sim$vova}%
\author[V.V. Uspenskij]{Vladimir V. Uspenskij}
\address{Department of Mathematics, 321 Morton Hall, Ohio University,
Athens, OH 45701, USA}
\email{uspensk@bing.math.ohiou.edu}
\date{Version as of November 21, 2000}

\keywords{Weakly almost periodic functions, precompact groups}
\begin{abstract} An action of a group $G$ on a compact space $X$ is
called weakly almost periodic if the orbit of every continuous function
on $X$ is weakly relatively compact in $C(X)$.
We observe that for a topological group $G$ the following are equivalent:
(i) every continuous action of $G$ on a compact space is
weakly almost periodic; (ii) $G$ is precompact.
For monothetic groups the result was previously obtained 
by Akin and Glasner, while for locally compact groups it has been
known for a long time.
\end{abstract}

\thanks{{\it 2000 Mathematical Subject Classification.} 
Primary 43A60, Secondary 22A05, 22C05, 22F05, 43A07.}

\maketitle

\section{Introduction}

Let a group $G$ act on a set $X$. Denote, as usual, by $l^\infty(X)$
the vector space of all bounded complex-valued functions on $X$
equipped with the supremum norm. 
A bounded function $f$ on $X$ is called 
{\it weakly almost periodic} (w.a.p. for short) if the
$G$-orbit of $f$ is weakly relatively compact in the
Banach space $l^\infty(X)$. 

Now suppose that a group $G$ acts by homeomorphisms on
a compact space $X$. The action is called {\it weakly almost
periodic} (or again {\it w.a.p.}) if every continuous function on
$X$ is weakly almost periodic. 

It follows easily from well-known results that if $G$ is a 
precompact topological group, then every continuous 
action of $G$ on a compact space is weakly almost periodic.
(Recall that a topological group is 
{\it precompact} if it is isomorphic to a subgroup of a compact group.)
Our main goal is to establish the converse.
\vskip .3cm

\noindent{\bf Main Theorem \ref{main}.} {\it 
For an arbitrary topological group $G$ the following
conditions are equivalent.
\begin{enumerate}
\item Every continuous action of $G$ on a compact space is
weakly almost periodic.
\item Every bounded right uniformly continuous function on $G$
is weakly almost periodic.
\item $G$ is precompact. \qed
\end{enumerate}}

In the case where $G$ is a {\it monothetic} group
(that is, contains an everywhere dense cyclic subgroup), the result
was established by Akin and Glasner \cite{AG}.

Our proof is independent and different in nature 
from that in \cite{AG}, and is
exploiting theory of invariant means and a result by Pachl \cite{Pachl}.

While surveying some basic notions and results from abstract topological
dynamics upon which our proof hinges, 
we give a new and simpler proof of the Ellis-Lawson
Joint Continuity Theorem (Theorem \ref{ellis} below).

\section{\label{weak}Weak almost periodicity and compactifications}

We begin with some background material about weak almost
periodicity, greatest ambits, and their 
significance in abstract topological dynamics. For more comprehensive 
information we refer the reader to e.g. \cite{Ru}. 

For a topological space $X$, let $C^b(X)$ stand for
the space of all bounded continuous
functions on $X$. Notice that $C^b(X)\sbs l^\infty(X)$ forms a
(weakly) closed vector subspace.
If $G$ acts on $X$ by homeomorphisms, then
the orbit of every function $f\in C^b(X)$ is contained in $C^b(X)$. 
It follows that $f$ is weakly almost periodic
if and only if its orbit is relatively
weakly compact in the Banach space $C^b(X)$.

In particular, considering the left and right actions of
a group $G$ on itself,
we can define left and right weakly almost periodic 
functions on $G$. These two notions are
actually equivalent \cite[Corollary 1.12]{Burckel}, so we can simply 
speak about w.a.p. functions on a group $G$. 

For a topological group $G$, denote by
$W(G)$ the space of all continuous w.a.p. functions on $G$. 
The space $W(G)$ is a commutative $C^\ast$-algebra and thus
is isomorphic to the algebra $C(G^w)$ of continuous functions on
a compact space $G^w$, the maximal ideal space of $W(G)$. 
Call $G^w$ the {\it weakly almost periodic (w.a.p.) 
compactification} of the topological group $G$.
There is a natural structure of a semitopological semigroup on $G^w$,
where `semitopological' means that the multiplication is 
separately continuous. 
The canonical semigroup homomorphism $G\to G^w$ 
is a universal object in the category
of continuous homomorphisms of $G$ to compact 
semitopological semigroups.

\begin{remark}
Under a compactification of a topological space $X$ we mean
a compact Hausdorff space $K$ together with a continuous
map $j:X\to K$ with a dense range. We do {\it not} require that $j$ be a
homeomorphic embedding. The long-standing problem of whether for every
topological group $G$ the canonical map $G\to G^w$ is a topological 
embedding has been
recently solved in the negative by the first named author. 
Namely, for the Polish group $G=\mathrm{Homeo}_+{\mathbb I}$
of all orientation-preserving homeomorphisms of the closed interval,
equipped with the compact-open topology, the  
w.a.p. compactification $G^w$
is a singleton (equivalently, every $f\in W(G)$ is constant) 
\cite{Me3}. 
\end{remark}

For a topological group $G$ the {\it greatest ambit} $\sS(G)$ is
the compactification of $G$ corresponding to the algebra $\RUC(G)$
of all bounded right uniformly continuous functions on $G$,
that is, the maximal ideal space of the latter $C^\ast$-algebra.
(A function $f\in C(G)$ is {\it right uniformly continuous} if
\[
\forall\e>0\,\exists V\in\sN(G)\,\forall x,y\in G\, (xy\obr\in V\implies
|f(y)-f(x)|<\e),
\]
where $\sN(G)$ is the filter of neighbourhoods of unity.)
There is a natural $G$-space structure on $\sS(G)$. (By a $G$-space
we mean a topological space $X$ equipped 
with a jointly continuous action 
$G\times X\to X$.) The canonical map $i:G\to\sS(G)$ is a homeomorphic 
embedding, and we will identify $G$ with $i(G)$. 
Let $e\in G\sbs \sS(G)$ be the
unity. The pair $(\sS(G),e)$ has the following universal property:
for every compact $G$-space $X$ with a distinguished point $p$ there exists
a unique $G$-map $\sS(G)\to X$ which sends $e$ to $p$. The multiplication
on $G$ extends to a multiplication on $\sS(G)$ such that every right shift
$r_a:\sS(G)\to \sS(G)$, $r_a(x)=xa$, is continuous. The shift $r_a$ is
the unique $G$-selfmap of $\sS(G)$ such that $r_a(e)=a$. 

Let $X$ be a Banach space. Denote
by $U(X)$ the semigroup of all linear operators on $X$ of norm $\le1$.
If $X$ is a reflexive space, then $U(S)$, 
equipped with the weak operator topology, 
becomes a compact semitopological semigroup. 
Every compact semitopological semigroup is isomorphic
to a closed subsemigroup of $U(X)$ for some reflexive Banach space $X$ 
\cite{Sh, Me3}.

Let again $X$ be a Banach space. 
For a subgroup $G$ of the group $\Is(X)$ of all linear  
isometries of $X$ denote by $G_w$ and $G_s$ the group $G$ equipped with
the weak and strong operator topology, respectively. For a wide class of 
Banach spaces $X$, including reflexive spaces, the two topologies
actually coincide: $G_w=G_s$ \cite{Me2}. 
If $X$ is a Hilbert space then $\Is(X)$ is 
weak-dense in $U(X)$ \cite{Us1}. 
Moreover, it follows from the results of \cite{Us1} that in this case 
$U(X)$ is the w.a.p. compactification of $\Is(X)$.

Each of the following three assertions easily implies the other two:
\begin{enumerate}
\item[]
(1) $W(G)\sbs\hbox{RUC}(G)$; 
\item[](2) there is a map $f:\sS(G)\to G^w$ such that
$j=fi$, where $i:G\to \sS(G)$ and $j:G\to G^w$ are the canonical maps;
\item[](3) the action of $G$ on $G^w$, induced by the semigroup homomorphism
$j:G\to G^w$, is jointly continuous.
\end{enumerate} 
All these are true for an arbitrary topological group $G$
(see, for example, \cite{Ru}), although this
is not at all obvious. Actually, the question of whether every w.a.p. 
function is uniformly continuous appeared 
on the top of the list of open problems in 
\cite{Burckel}. The assertion (3) above follows from a more general
result of J. Lawson, which in itself is a corollary of a 
fundamental {\it Ellis-Lawson 
Joint Continuity Theorem} \cite{La}. We include here a new ``soft'' proof
of this result.

\begin{thm}[\cite{La}, and \cite{Ellis} for $S=G$] 
Let $S$ be a compact semitopological semigroup with unity. Let $G$ be the
group of invertible elements in $S$. Then $G$ with the induced topology
is a topological group, and the action of $G$ on $S$ by translations is
jointly continuous.
\label{ellis}
\end{thm}

\begin{proof}
As we have noted before, there exists a reflexive Banach space $X$ 
such that $S$ is a subsemigroup of $U(X)$. We can therefore
suppose that $S=U(X), G=\Is(X)$. We have $G_s=G_w$ \cite{Me2}. 
Since $G_s$ is a topological group, so is  
$G_w$, that is, $G$ with the topology it inherits as a subspace of $S$.
To see that the action of $G$ on $U(X)$ is (jointly) continuous, it suffices
to prove that the action of $G=\Is(X)_s$ on the dual Banach 
space $X^*$ is continuous.
This fact has been proved in \cite[Corollary 6.9]{Me0}, and
easily follows also by observing that the topological groups 
$\Is_w(X)=\Is_s(X)$ and 
$\Is_w(X^*)=\Is_s(X^*)$ are naturally isomorphic.
\end{proof}

\section{Invariant means}

Let $G$ be an arbitrary group and let $A$ be
a $C^\ast$-subalgebra of $l^\infty(G)$. A {\it mean} on $A$ is
a complex-valued linear
functional $\phi$ on $A$ that is positive
(that is, $\phi(f)\geq 0$ whenever $f\geq 0$) and takes
the constant function $1$ to $1$. (In $C^\ast$-algebra theory,
one uses the term {\it state} instead.) Every mean is
automatically bounded of norm one
(cf. \cite{Sa}, Prop. 1.5.1).
A mean $\phi$ is {\it left-invariant} if $\phi(_gf)=\phi(f)$ for
every $f\in A$ and every $g\in G$,
where $_gf$ denotes the left translate of the function $f$ by $g$.
In a similar way one defines {\it right-invariant means} on $A$. 
If a mean $\phi$ is both left and right invariant, we will refer to
it as a {\it bi-invariant} mean.

For instance, there is always a bi-invariant mean on the algebra
$W(G)$, as the following well-known result asserts.
(See e.g.  \cite[Corollary 1.26]{Burckel}.)

\begin{thm}[Ryll-Nardzewski] For every topological group $G$,
there is a unique bi-invariant mean on the algebra $W(G)$. 
Moreover, such a mean is the unique left-invariant mean on
$W(G)$ as well.
\label{ryll}
\qed
\end{thm}

\begin{remark}
The cited corollary actually only
asserts that there exists a unique two-sided invariant mean on $W(G)$.
To show the uniqueness in the class of left-invariant means, note that
for every $f\in W(G)$ there is a constant that can be uniformly
approximated by convex combinations of left translates of $f$ 
\cite[Theorem~1.25 and Corollary~1.26]{Burckel}. The value of any 
left-invariant mean at $f$ must be equal to such a constant, whence
the uniqueness.
\end{remark}

On the contrary, the larger algebra $\RUC(G)$ need not
in general support even a left-invariant mean.
If there is such a mean on $\RUC(G)$, the topological group
$G$ is called {\it amenable}.  
Equivalently, $G$ is amenable if for every
compact convex set $K$ (lying in some locally convex space) and every 
continuous action of $G$ on $K$ by affine maps 
there exists a $G$-fixed point in $K$. See e.g. \cite{Pa}.

Left-invariant means on $\RUC(G)$ can be identified with
left-invariant probability measures on the greatest ambit 
$\sS(G)$. The collection of all such left-invariant means, 
equipped
with the weak$^\ast$ topology, forms a compact space, which we will
denote by $\hbox{LIM}_G$. 

If $\hbox{LIM}_G$ consists of a single point, we
will call the topological group $G$
{\it uniquely amenable.} Every compact group provides an 
obvious example of a uniquely amenable group, for which the
unique invariant mean comes from the Haar measure. 
It is also obvious that every precompact group $G$ is uniquely
amenable as well, since the algebras $\RUC(G)$ and
$\RUC(\widehat G)$ are canonically isomorphic, where $\widehat G$
denotes the compact completion of $G$. This observation can be,
at least partially, inverted, as the following result shows.

\begin{thm}[Pachl, \cite{Pachl}] 
Let $G$ be a separable metrizable group.
If the compactum $\hbox{LIM}_G$ contains
a $G_\delta$-point, then $G$ is a precompact group. \qed
\label{pach}
\end{thm}

\begin{corol}
Every uniquely amenable separable metrizable topological group
is precompact. \qed
\label{uniquely}
\end{corol}

It remains unclear if the above Corollary remains true for
arbitrary topological groups. 

\begin{question}
Is it true that every uniquely amenable topological 
group is precompact?
\label{conj}
\end{question}

If the answer is in the affirmative, then our Main Theorem
\ref{main} below follows at once, and
Proposition~\ref{invar} and Lemma~\ref{lemma} are not needed.
We can only provide a positive answer to Question \ref{conj}
under the additional assumption that $G$ is $\o$-bounded. 

A topological group
is called $\o${\it -bounded} if for every $U\in \sN(G)$ 
there exists a countable $A\sbs G$ such
that $AU=G$, or, equivalently, if $G$ is 
isomorphic to a subgroup of the product
of separable metrizable groups \cite{Arh, Gu}. 

\begin{lemma}
The property of being uniquely amenable is preserved 
by continuous homomorphisms.
\end{lemma}

\begin{proof}
If $G$ is amenable and $f:K\to L$ is an affine onto $G$-map between 
compact convex $G$-spaces, then for every $G$-fixed 
point $y\in L$ there exists
a $G$-fixed point $x\in K$ such that $f(x)=y$. 
Indeed, the set $f\obr(y)$ is a
compact convex $G$-space and hence contains a $G$-fixed point.

Let $f:G\to H$ be an onto homomorphism. Assume that $G$ is uniquely amenable.
Applying the remark of the previous paragraph 
to the sets of probability measures
on $\sS(G)$ and $\sS(H)$, we see that every 
left-invariant measure on $\sS(H)$
is the image under $f$ of a left-invariant 
measure on $\sS(G)$. Since the latter
is unique, it follows that the left-invariant measure on $\sS(H)$ 
is also unique.
\end{proof}

Combining this result with Theorem \ref{pach}, we obtain:

\begin{corol}
Every uniquely amenable $\o$-bounded group is precompact. \qed
\end{corol} 

\section{The main result}

As we have seen in Section \ref{weak}, 
$W(G)\sbs \hbox{RUC}^b(G)$ for every topological group $G$. If $G$
is compact, then $W(G)= \hbox{RUC}^b(G)$. Indeed, it is well known
that every continuous function $f$ on a compact group is
{\it almost periodic,} that is, the orbit of $f$ under left (or
right) translations is even {\it norm relatively compact} in $C(G)$. 
The same is true of every precompact group $G$. 

The aim of this paper is to prove that the equality 
$W(G)= \hbox{RUC}^b(G)$ is
actually {\it equivalent} to $G$ being precompact. 
Such fact is well known for locally compact groups; see 
\cite{Dz} and \cite{Ch} (the latter for the case of 
locally compact groups $G$ having small invariant neighbourhoods).
For monothetic groups (that is, those --- not necessarily
locally compact --- topological groups
containing an everywhere dense cyclic subgroup) 
this result was recently established by Akin and Glasner, 
whose preprint \cite{AG} has stimulated our present research.

We begin with the following auxiliary (and well known) statement.

\begin{prop}
For every topological group $G$ the following are equivalent:
\begin{enumerate}
\item $W(G)=\hbox{RUC}^b(G)$, that is, every right uniformly continuous
function on $G$ is weakly almost periodic;
\item $\sS(G)=G^w$;
\item the natural multiplication on $\sS(G)$ is separately continuous (on both
sides); 
\item if $X$ is a compact $G$-space and $f\in C(X)$, then $f$ is weakly 
almost periodic.
\end{enumerate}
\label{equiv}
\end{prop}

\begin{proof}
The preceding discussion shows that the conditions (1) -- (3) are equivalent.
Applying (4) to $X=\sS(G)$, we see that $(4)\implies(1)$. The implication
$(1)\implies(4)$ can be established directly. We omit the proof, since we
show below (Theorem~\ref{main})
that (1) actually implies that $G$ is precompact, and for 
precompact groups (4) is straightforward: just note that the $G$-space
structure on a compact space $X$ extends to a 
$\widehat G$-space structure.
Of course one has to keep in mind that in our proofs the
implication $(1)\implies(4)$ is never used.
\end{proof}

\begin{defin} Following \cite{AG},
we will say provisionally that a topological group 
$G$ is {\it weakly almost periodic} (w.a.p.) 
if it satisfies the equivalent conditions of Proposition~\ref{equiv}.
\end{defin}

\begin{lemma}
The class of w.a.p. groups is closed under forming continuous
homomorphic images and topological subgroups.
\label{invar}
\end{lemma}

\begin{proof}
Let $G$ be a w.a.p. topological group, 
and let $h:G\to H$ be an onto homomorphism. If 
$f\in\hbox{RUC}^b(H)$, then $fh\in \hbox{RUC}^b(G)=W(G)$. The map
$h^*:C(H)\to C(G)$ is an isometric  embedding for the norm topologies
and hence a homeomorphic embedding for the weak topologies. Since $h^*$
sends the orbit of $f$ to the orbit of $fh$ and the orbit of $fh$ is
weakly relatively compact, so is the orbit of $f$. Thus RUC$^b(H)\sbs W(H)$.

Now let $H$ be a subgroup of $G$ and $f\in\hbox{RUC}^b(H)$. A uniformly
continuous bounded function defined on a subspace of a uniform space can
be extended to a uniformly continuous bounded function defined on the entire
space \cite[8.5.6]{Eng},
\cite{Ka}. Since the right uniformity of $H$ is induced by the right
uniformity of $G$, there exists a function $g\in\hbox{RUC}^b(G)$ such that 
$g\restriction H=f$. By assumption, $g\in W(G)$. 
The restriction map $C(G)\to C(H)$ sends weakly relatively compact
subsets of $C(G)$ to weakly relatively compact subsets of $C(H)$,
while the image of the $G$-orbit of $g$ contains the $H$-orbit of $f$.
We conclude that $f\in W(H)$.
\end{proof}

\begin{lemma}
Let $\sP$ be a class of topological groups closed under subgroups and 
homomorphic images. If not all groups in $\sP$ are precompact, then
there exists a countable metrizable group in $\sP$ which is not precompact.
\label{lemma}
\end{lemma}

\begin{proof}
Let $G\in \sP$ be non-precompact. There exists an infinite set
$A\sbs G$ and a neighbourhood $U$ of unity such that 
the family $\{Ua:a\in A\}$
is disjoint. The subgroup of $G$ generated by $A$ is countable and 
non-precompact. Thus we may assume that $G$ is countable. 
Every countable group, being $\o$-bounded,
is isomorphic to a subgroup of the product $\prod G_\a$ 
of countable metrizable groups.
We may assume that the projections of $G$ to the factors $G_\a$ are onto
and hence all the $G_\a$'s are in $\sP$.
Since the class of precompact groups is closed 
under products and subgroups, at least one
of the factors $G_\a$'s has to be non-precompact. 
\end{proof}

We now turn to the main result of the paper.

\begin{thm} For an arbitrary topological group $G$ the following
conditions are equivalent.
\begin{enumerate}
\item Every continuous action of $G$ on a compact space is
weakly almost periodic.
\item Every bounded right uniformly continuous function on $G$
is weakly almost periodic.
\item $G$ is precompact.
\end{enumerate}
\label{main}
\end{thm}

\begin{proof}
Only (2)$\implies$(3) requires a proof. Let $W(G)=\RUC(G)$. 
As an immediate consequence of this assumption and 
Ryll-Nardzewski Theorem \ref{ryll}, the topological group $G$
is uniquely amenable. 

Assume first that $G$ is separable metrizable. 
According to Pachl's Corollary \ref{uniquely},
$G$ is precompact.  

The general case is now being reduced to the case of a separable metrizable
group with the aid of Proposition~\ref{invar} and Lemma~\ref{lemma}.
\end{proof}

\subsection*{Acknowledgements} 
The second named author (V.P.) is grateful to the Department of
Mathematics of Tel-Aviv University for hospitality extended
during his visit in April 1999, and to Eli Glasner for useful
discussions, in particular for asking the question on the
validity of what became now Theorem \ref{main} of this note,
and later forwarding a copy of the preprint \cite{AG}.

\bibliographystyle{amsplain}

\begin{thebibliography}{10}


\bibitem{AG} E. Akin and E. Glasner, {\it Residual properties and
almost equicontinuity,} Tel-Aviv University preprint, October 8, 2000.

\bibitem{Arh} A.V. Arhangel'ski\u\i, 
{\it Classes of topological groups,} Russ. Math. Surveys
{\bf 36}:3 (1981), 151--174.

\bibitem{Burckel} R.B. Burckel,
{\it Weakly almost periodic functions on semigroups,}
Gordon and Breach Science Publishers, 
New York-London-Paris, 1970, ix+118 pp.

\bibitem{Ch} 
C. Chou, {\it Weakly almost periodic functions and almost 
convergent functions on a group}, Trans. Amer. Math. Soc.
{\bf 206} (1975), 175-200.

\bibitem{Dz} 
H.A.M. Dzinotyiweyi, {\it Nonseparability of quotient 
spaces of function algebras on topological semigroups}, 
Trans. of the Amer. Math. Soc., {\bf 272}:1 (1982), 223-235.

\bibitem{Ellis} R. Ellis, {\it Locally compact transformation groups}, 
Duke Math. J. vol. {\bf 24}, (1957), 119-125. 

\bibitem{Eng} R. Engelking, {\it General Topology,}
PWN, Warczawa, 1977.

\bibitem{Gu} I. Guran, 
{\it On topological groups close to being Lindelof}, 
Soviet Math. Dokl., {\bf 23} (1981), 173-175.

\bibitem{Ka} M. Katetov, {\it On 
real-valued functions in topological spaces}, 
Fund. Math. {\bf 38}, 85-91, and correction {\bf 40} (1953), 203-205. 
 
\bibitem{La} J.D. Lawson, 
{\it Joint continuity in semitopological semigroups}, 
Illinois J. Math., {\bf 18} (1974), 275--285.

\bibitem{Me0} M.G. Megrelishvili, {\it Fragmentability and continuity of 
semigroup actions}, Semigroup Forum, {\bf 57} (1998), 101-126. 

\bibitem{Me1} M.G. Megrelishvili, 
{\it Eberlein groups and compact semitopological semigroups}, 
Bar-Ilan University preprint, 1998.

\bibitem{Me2} M.G. Megrelishvili, 
{\it Operator topologies and reflexive 
representability of groups}, 
BIUMCS 98/52, Bar-Ilan University preprint, 1998. 

\bibitem{Me3} M.G. Megrelishvili, 
{\it Every semitopological semigroup compactification of the group 
$H_+[0,1]$ is trivial},
Semigroup Forum (to appear).

\bibitem{Pachl} J.K. Pachl, {\it Uniform measures on topological groups,}
Compositio Math. {\bf 45} (1982), 385--392. 
 
\bibitem{Pa} A.T. Paterson, 
{\it Amenability,} Math. Surveys and Monographs
{\bf 29}, Amer. Math. Soc., Providence, RI, 1988.

\bibitem{Ru} W. Ruppert, 
{\it Compact semitopological semigroups: An intrinsic theory}, 
Lecture Notes in Math., {\bf 1079}, Springer-Verlag, 1984.

\bibitem{Sa} S. Sakai, {\it $C^\ast$-Algebras and $W^\ast$-Algebras,}
Springer-Verlag, Berlin--Neidelberg--NY, 1971; Reprinted,
Springer, 1998.

\bibitem{Sh} A. Shtern, {\it Compact semitopological semigroups and 
reflexive representability of topological groups}, 
Russian J. of Math. Physics, 
{\bf 2} (1994), 131-132.

\bibitem{Us1} V.V. Uspenskij, 
{\it The Roelcke compactification of unitary groups}, 
in: Abelian groups, module theory, and topology, 
Proceedings in honor of Adalberto Orsatti's 60th birthday 
(D.  Dikranjan, L. Salce, eds.), 
Lecture notes in pure and applied mathematics, 
Marcel Dekker, New York e.a., {\bf 201} (1998), 411-419.

\end{thebibliography}

\end{document}